\documentclass[10pt]{amsart}

\usepackage{amsfonts}
\usepackage[]{amsmath}
\usepackage[]{amssymb}
\usepackage[]{amsthm}
\usepackage[]{amsfonts}

\usepackage{setspace}
\newtheorem  {theorem}            {Theorem}
\newtheorem  {lemma}              {Lemma}
\newtheorem  {definition}         {Definition}
\newtheorem  {corollary}          {Corollary}

\newtheorem* {conjecture*}        {Conjecture}
\newtheorem* {theorem*}           {Theorem}
\newtheorem* {acknowledgements*}  {Acknowledgements}

\newcommand {\Ric}  {\operatorname{Ric}}
\newcommand {\R}    {\operatorname{R}}
\newcommand {\A}    {\operatorname{A}}
\newcommand {\mc}   {\operatorname{H}}
\newcommand {\area} {\operatorname{area}}
\newcommand {\vol}  {\operatorname{vol}}
\newcommand {\supp} {\operatorname{supp}}

\title[Isoperimetric surfaces in 3-manifolds] {The size of isoperimetric
surfaces in $3$-manifolds and a rigidity result for the upper hemisphere}

\author{Michael Eichmair}
\address{Michael Eichmair, Department of Mathematics, MIT, 77 Massachusetts
Avenue, Cambridge MA 02139-4307, USA }
\email{eichmair@math.mit.edu}

\begin{document}

\maketitle

\begin{abstract} We characterize the standard $\mathbb{S}^3$ as the closed
Ricci-positive $3$-manifold with scalar curvature at least $6$ having
isoperimetric surfaces of largest area: $4\pi$. As a corollary we answer in the
affirmative an interesting special case of a conjecture of M. Min-Oo's on the
scalar curvature rigidity of the upper hemisphere.
\end{abstract}


\section{Introduction}

The following rigidity result for the unit ball of $(\mathbb{R}^3, \delta)$ is
a well-known consequence of the positive mass theorem (see \cite{PMTI},
\cite{Witten}, \cite{Miao}, \cite{ShiTam}, and \cite{HangWang}):

\begin{theorem*} Let $(M^3, g)$ be a compact orientable Riemannian $3$-manifold
with non-negative scalar curvature and boundary isometric to round
$\mathbb{S}^{2}$ with mean curvature equal to $2$. Then $(M^3, g)$ is isometric
to the unit ball of $(\mathbb{R}^3, \delta)$. \end{theorem*}

\noindent The theorem asserts that there are no compact deformations of the
Euclidean metric within the class of non-negative scalar curvature metrics. The
analogous rigidity statement for hyperbolic space was proven in \cite{Min-Oo},
\cite{Wang}, \cite{ChruscielHerzlich}, \cite{AnderssonCaiGalloway} by
establishing appropriate versions of the positive mass theorem in this context.
In \cite{Min-Oo98} M. Min-Oo raised the following question:

\begin{conjecture*} [Min-Oo] Let $(M^n, g)$ be an $n$-dimensional compact
orientable Riemannian manifold with scalar curvature $\R\geq n(n-1)$ and
totally geodesic boundary isometric to round $\mathbb{S}^{n-1}$. Then $(M^n,
g)$ is isometric to the round hemisphere $\mathbb{S}^n_+$.
\end{conjecture*}

\noindent We refer the reader to \cite{Min-Oo98}, \cite{HangWang}, and
\cite{HangWang2} for more background and context for Min-Oo's conjecture.
Recall that the standard metrics on $\mathbb{R}^n, \mathbb{H}^n,$ and
$\mathbb{S}^n$ are all static (the linearization of the scalar curvature map
about the standard metrics on these spaces has non-trivial cokernel). By the
work of J. Corvino (Theorem 4 in \cite{Corvino}), one can always locally deform
the scalar curvature in \emph{any} direction if the underlying metric is not
static:

\begin{theorem*} [Corvino] Let $\Omega$ be a smooth domain compactly contained
in a Riemannian manifold $(M^n, g_0)$. Suppose that the linearization $L_{g_0}$
of the scalar curvature map $\R : \mathcal{C}^\infty(\Omega) \to
\mathcal{C}^\infty(\Omega)$ at $g_0$ has an injective formal
$\mathcal{L}^2$-adjoint $L^*_{g_0}$, where we consider $L_{g_0} :
H^2_{\text{loc}}(\Omega) \to \mathcal{L}^2_{\text{loc}}(\Omega)$. Then for
every smooth function $S$ on $\Omega$ sufficiently close to the scalar
curvature $\R(g_0)$ and equal to $\R(g_0)$ near $\partial \Omega$, there exists
a smooth metric $g$ on $\Omega$ such that $\R(g) = S$ and so that $g\equiv g_0$
outside $\Omega$. \end{theorem*}

\noindent As was noted in \cite{HangWang}, by Corvino's theorem, staticity
appears as an obstruction to finding potential counterexamples to Min-Oo's
conjecture near the round metric on $\mathbb{S}^n_+$. However, in \cite
{HangWang}, F. Hang and X. Wang show that if one moves from the hemisphere
$\mathbb{S}^n_+$ to a larger geodesic ball of $\mathbb{S}^n$, there even are
metrics conformally related to the round metric with scalar curvature $\geq
n(n-1)$ and standard boundary geometry:

\begin{theorem*} [Hang and Wang] For any $r \in (\frac{\pi}{2}, \pi)$ there is a
smooth metric $g = e^{2 \phi} g_{\mathbb{S}^3}$ on $\mathbb{S}^3$ such that (a)
$\R_g \geq 6$, (b) $\phi$ is not identically $0$, and (c) $\supp(\phi) \subset
B(N, r)$, where $N$ is a fixed point in $\mathbb{S}^3$.
\end{theorem*}

\noindent By contrast, in the same paper the authors establish Min-Oo's
conjecture among conformal deformations:

\begin{theorem*} [Hang and Wang] Let $g = e^{2 \varphi} g_{\mathbb{S}^n}$ be a
$\mathcal{C}^2$-metric on the upper hemisphere $\mathbb{S}^n_+$ satisfying the
assumptions (a) $\R_g \geq n(n-1)$ and (b) the boundary is totally geodesic and
isometric to standard $\mathbb{S}^{n-1}$. Then $g$ is isometric to
$g_{\mathbb{S}^n}$.
\end{theorem*}

In a recent paper \cite{HangWang2}, F. Hang and X. Wang use Raleigh's B\^ochner-type formula on manifolds with boundary to show that Min-Oo's conjecture holds true in all dimensions if one adds $\Ric_M \geq (n-1)g$ to the hypotheses. \\

The main result of this work answers Min-Oo's conjecture in the affirmative in
the case where $n=3$, $\Ric_M > 0$, and the boundary is an isoperimetric
surface (Theorem \ref{thm: Min-Oo}):

\begin{theorem*} Let $(M^3, g)$ be a compact orientable Riemannian $3$-manifold
with scalar curvature $\R_M \geq 6$, Ricci curvature $\Ric_M > 0$ and totally
geodesic boundary $\partial M^3$. If $area (\partial M^3) \geq 4 \pi$ and
$\partial M^3$ is an isoperimetric surface for the doubled manifold $(\hat M^3,
\hat{g})$, then $(M^3, g)$ is isometric to the upper hemisphere $\mathbb{S}^3_+$. \end {theorem*}

\noindent We point out that it is not necessary here to assume that $\partial M^3$ is round. In fact, the author is not aware of counterexamples to Min-Oo's conjecture when the original condition on the inner geometry of $\partial M$ is weakened to a lower bound on its area, as above. \\

\noindent Our main contribution here, which quickly leads to a proof of the above
theorem, is to characterize round $\mathbb{S}^3$ as the unique Ricci-positive
$3$-manifold with scalar curvature at least $6$ that admits isoperimetric
surfaces of largest area (Theorem \ref{thm: conclusion}):

\begin{theorem*} Let $(M^3, g)$ be a closed Riemannian manifold with $\R_M \geq
6$ and $\Ric_M > 0$. Then all isoperimetric surfaces of $M^3$ have area
strictly less than $4\pi$ unless $(M^3, g)$ is isometric to $\mathbb{S}^3$.
\end {theorem*}

\noindent There are two steps in the proof of this theorem. The crucial first ingredient is the monotonicity of a certain isoperimetric mass that H. Bray discovered in his thesis \cite{HughBrayThesis}. H. Bray's arguments imply directly that the isoperimetric profile of $(M^3, g)$ coincides with that of round $\mathbb{S}^3$. We recall that the isoperimetric surfaces of $\mathbb{S}^3$ are the geodesic balls (this follows from a symmetrization argument, see for example \cite{Ros}). Second, to show that $(M^3, g)$ is in fact round itself, we employ a delicate comparison with the isoperimetric ratio of small geodesic balls in $M^3$, using the Taylor expansion in \cite{GrayVanhecke} for their volume. \\

\noindent A curious feature of this proof is how the global information in the isoperimetric assumption implies complete rigidity of the local geometry. H. Bray's results in \cite{HughBrayThesis} have not been published. In the next section, we will summarize the required portion of his method to explain how it applies in this paper. \\

Note that the upper bound of $4\pi$ on the size of isoperimetric surfaces in the preceding theorem is already contained in the work of D. Christodoulou and S.-T. Yau \cite{ChristodoulouYau} and comes out of a Hersch-type choice of test functions in the stability inequality:

\begin{theorem*} [Christodoulou and Yau] \label{thm: ChristodoulouYau} Let $(M^3,
g)$ be a Riemannian $3$-manifold and let $\Sigma^2 \subset M^3$ be an immersed
closed (weakly) stable constant mean curvature surface of genus $0$. Then
\begin{equation*} 16 \pi \geq  \int_{\Sigma} \mc_{\Sigma}^2 + \frac
{2}{3}\int_{\Sigma} \R_M,
\end{equation*} where $\R_M$ is the scalar curvature of $M^3$.
\end{theorem*}

We point out that highly successful notions of mass and quasi-local mass, which are derived from the isoperimetric defect (from Schwarzschild) of outward minimizing surfaces in asymptotically flat  initial data sets with non-negative scalar curvature, have been introduced and studied by G. Huisken, who has used them to prove positive mass and Penrose-type theorems, including the rigidity cases. His work is largely unpublished at this point, but see the report \cite{Huisken}.  

\begin{acknowledgements*} This paper forms part of my thesis at Stanford University. I
would like to sincerely thank my advisor Richard Schoen for his constant
support and encouragement and for suggesting this problem to me. I am very
grateful to Hubert Bray for his interest in this work and for many stimulating and fun
discussions.
\end{acknowledgements*}


\section{Proofs}

We review some standard results regarding the isoperimetric profile function of a closed Riemannian manifold, starting with its definition. We refer to the excellent survey article \cite{Ros} and the paper \cite{BavardPansu} for the history and basic properties of the isoperimetric profile, as well as for further references on this topic.

\begin{definition} [\cite{BavardPansu}] Given a closed Riemannian $3$-manifold $(M^3, g)$, define its
isoperimetric profile function $I :[0, \vol(M^3)] \to \mathbb{R}$ by
\begin{displaymath} I(V) = \inf \{ \area (\partial \Omega) : \Omega \subseteq
M^3 \textrm{ region with } \\ \vol(\Omega) = V \}. \end{displaymath}
\end{definition}

\noindent It is a classical result in geometric measure theory that the infimum in the definition of the isoperimetric profile is achieved by ``isoperimetric regions" $\Omega \subset M^3$ with smooth embedded boundaries $\Sigma^2 = \partial \Omega$. Surfaces $\Sigma^2 \subset M^3$ arising in this way are called isoperimetric surfaces. A standard first variation argument yields that isoperimetric surfaces have constant mean curvature (the same constant for all connected components). The definition implies that $I(V)$ is symmetric with respect to $\frac{1}{2} \vol(M^3)$. \\

The following basic regularity of the isoperimetric profile was established in \cite{BavardPansu}:

\begin{lemma} [\cite{BavardPansu}] \label{lem: BavardPansu} Given $(M^3, g)$ closed and $V \in (0, \vol(M^3))$, let $\Omega \subset M^3$ be an isoperimetric region with $\vol(\Omega) = V$ and denote $\partial \Omega = \Sigma^2$. Write $\A_\Sigma$, $\mc_{\Sigma}$ for the second fundamental form and (constant) mean curvature of $\Sigma^2$ computed with respect to its outward normal $\vec\nu$. The isoperimetric profile has the following regularity:

a) $I$ has left and right derivatives at $V$, and $I'^+(V) \leq \mc_\Sigma \leq
I'^-(V)$.

b) $I''(V) I(V)^2 + \int_\Sigma (\Ric_M(\vec \nu,\vec\nu) + |\A_\Sigma|^2) \leq
0$ holds in the sense of comparison functions.
\end{lemma}

\noindent The lemma asserts that for every $V\in(0, \vol(M^3))$ there exists a smooth function $I_{V}$ defined in a neighbourhood of $V$ such that $I_{V}(V) = I(V)$, $I_{V} \geq I$, and $I_{V}''(V)I_{V}^2(V) + \int_\Sigma (\Ric_M(\vec\nu,\vec\nu) + |\A_\Sigma|^2) \leq 0$. A well-known and immediate consequence of part b) is that the isoperimetric profile is concave when $(M^3, g)$ has non-negative Ricci curvature. \\

In his thesis, H. Bray proved a scalar curvature based volume comparison theorem for $3$-manifolds. More precisely, he showed that if a closed $3$-manifold $(M^3, g)$ has scalar curvature bounded below by $6$ and Ricci curvature bounded below by $\varepsilon (2  g)$ for some $\varepsilon \in (0, 1)$, then its volume is bounded above by the volume of the round unit sphere $\mathbb{S}^3$ times a constant $\alpha(\varepsilon)$, where $\alpha(\varepsilon) = 1$ for large enough $\varepsilon \in (0, 1)$. The techniques from \cite{HughBrayThesis}, which are crucial in our proof of Theorem \ref{thm: conclusion}, have not been published. For convenient reference, and in order to explain how the arguments from \cite{HughBrayThesis} apply directly in our context, we summarize several statements and proofs from H. Bray's thesis below: 

\begin{definition} [\cite{HughBrayThesis}] Let $(M^3, g)$ be a closed Riemannian $3$-manifold with
scalar curvature $\R_M \geq 6$. The adapted Hawking $m_H : (0, \vol(M^3)) \to
\mathbb{R}$ is defined in terms of the isoperimetric profile by \[m_H(V) =
\sqrt{I(V)} \Big(16\pi - 4I(V) - I(V)I'^+(V)^2 \Big).\]
\end{definition}

It was proven in \cite{HughBrayThesis} that isoperimetric surfaces are connected if the ambient manifold has positive Ricci curvature. In conjunction with the Gauss-Bonnet theorem and the estimate in part b) of Lemma \ref{lem: BavardPansu}, H. Bray obtained the following result:

\begin{lemma} [\cite{HughBrayThesis}] \label{lem: basic estimate} Assume that $(M^3, g)$ has positive Ricci curvature. Then \[ I''(V) \leq I''_V(V) = - \frac {\int_\Sigma |\A_\Sigma|^2 + \Ric_{M}(\vec\nu, \vec\nu)}{I(V)^2} \leq \frac {4\pi}{I(V)^2} - \frac {3 I'^+(V)^2}{4I(V)} - \frac {\int_\Sigma \R_M}{2I(V)^2}\] for every isoperimetric surface $\Sigma^2$ corresponding to the volume $V$. \end{lemma}

The monotonicity of the adapted Hawking mass $m_H$ in \cite{HughBrayThesis} is crucial:

\begin{lemma} [\cite{HughBrayThesis}] \label{lem: monotonicity} Let $(M^3, g)$ be a closed
Riemannian $3$-manifold with $\R_M \geq 6$. Then, as a distribution, $m_H'\geq0$
on any connected subinterval of $(0, \vol(M^3)/2))$ on which (a) every volume
is realized by some isoperimetric region with connected boundary, and (b) the
isoperimetric profile $I$ is nondecreasing. In particular, if $\Ric_M > 0$, then $m_H$ is nondecreasing on the interval $(0, \frac{1}{2} \vol(M^3))$.

\begin{proof}
For $\delta \neq 0$ define the difference quotient operator
$\Delta_{\delta}|_{V}f = \delta^{-1}(f(V+\delta) - f(V))$. Recall that
$\Delta_\delta$ obeys a product rule and has formal adjoint
$-\Delta_{-\delta}$. Moreover, if $f\geq g$ in a neighbourhood of a point $V$
with $f(V) = g(V)$, then $\Delta_{-\delta}|_{V} \Delta_{\delta} f \geq
\Delta_{-\delta}|_{V} \Delta_{\delta} g$ for $\delta \neq 0$ sufficiently
small. \\

\noindent Let $ 0 \leq \phi \in C_c^1 (a, b)$ be a nonnegative test function
compactly supported in an interval $(a, b) \subset (0, \vol(M^3))$ as in the
hypothesis of the lemma. One computes
\begin{eqnarray*} - \int \phi' m_H & = & - \lim_{\delta \to 0} \int
(\Delta_\delta \phi) \sqrt{I}\big(16 \pi - 4I - (\Delta_\delta I)^2 I \big)
 \\ & = & \lim_{\delta \to 0} \int  \phi \Delta_{-\delta} \Big(\sqrt{I}\big(16
\pi - 4I - (\Delta_\delta I)^2 I \big)\Big) \\ & = & \int 2 \phi I' I^{\frac 3
2} \big( \frac {4\pi}{I^2} - \frac {3 I'^2}{4I} - \frac{3}{I}\big) +
\lim_{\delta \to 0} \int 2\phi I' I^{\frac 3 2}
(-\Delta_{-\delta}(\Delta_\delta I)) \\ & \geq & \int 2 \phi I' I^{\frac 3 2}
\Big( \frac {4\pi}{I^2} - \frac {3 I'^2}{4I} - \frac{3}{I} - \limsup_{\delta
\to 0} \Delta_{-\delta}(\Delta_\delta I)\Big). \end{eqnarray*}

\noindent Here, it was used that $I'^+(V) = I'^-(V)$, except possibly at
countably many values of $V$, and the fact that $\Delta_{-\delta}(\Delta_\delta
I) \leq 0$, which is just the concavity of $I$. By Lemma \ref{lem: basic estimate}, \begin{eqnarray*}
\limsup_{\delta \to 0} \Delta_{-\delta}|_V(\Delta_\delta I) \leq I''_V(V) &
\leq & \frac {4\pi}{I^2(V)} - \frac {3 I'^+(V)^2}{4I(V)} -
\frac{\int_{{\Sigma}_V} \R_M} {2 I(V)^2} \\ & \leq & \frac {4\pi}{I^2(V)} -
\frac {3 I'^+(V)^2}{4I(V)} - \frac{3} {I(V)} \end{eqnarray*} for almost every
$V$ since again $I'^+(V) = I'^-(V)$ for all but countably many values of $V$. Together with the hypothesis that $I'^+ \geq 0$, this implies that
$-\int \phi' m_H \geq 0$. \\

By the concavity of $I$ when $\Ric_M > 0$ and its symmetry with respect to $\frac{1}{2} \vol(M^3)$, $I_M$ is increasing on $(0, \vol(M^3)/2)$, hence, by the above, so is $m_H$. 
\end{proof}
\end{lemma}

The following corollary characterizes equality in Lemma \ref{lem: monotonicity}. Its proof is implicit in the proof of the scalar curvature rigidity theorem in \cite{HughBrayThesis}. 

\begin{corollary} [\cite{HughBrayThesis}] \label{cor: isoperimetric profiles agree} Let $(M^3, g)$ have $\R_M \geq 6$ and $\Ric_M > 0$. If a volume $V_0 \in (0, \vol(M^3)/2)]$ is such that $m_H(V_0) \leq 0$, then the isoperimetric profile of $M^3$ coincides with that of $\mathbb{S}^3$ on $[0, V_0]$.

\begin{proof} Lemma \ref{lem: monotonicity} implies that $m_H' \geq 0$ as a distribution on $(0, \vol(M^3)/2)$, and hence is a non-decreasing function at least where $m_H$ is continuous. Since $I$ is concave, $I'^+$ can only jump down and it follows that $m_H$ is in fact non-decreasing on all of $(0, \vol(M^3)/2]$. Moreover, $\lim_{V \to 0 } m_H (V) = 0$ (one can use the result of Christodoulou and Yau in the introduction to argue this) and hence if $m_H\big(V_0) = 0$, then $m_H$ has to vanish identically on $(0, V_0]$. This means that $16\pi - 4I(V) - I(V)I'^+(V)^2 \equiv 0$ on this interval. From the continuity of $I$ it follows that $I' \equiv I'^+$ is continuous and that $I$ is a classical solution of the ODE \[ I' = \sqrt {\frac {16\pi - 4I}{I}}\] on $(0, V_0]$. Since $\Ric_M >0 $, $I$ is strictly concave. It is easy to see now that there exists a unique solution with $I(0)=0$ (introducing the function $r = r(V)$ defined implicitly by $I(V) = 4\pi\sin^2(r)$ helps integrate the separated equation where $I(V) < 4 \pi$). An easy computation shows that the adapted Hawking mass of $\mathbb{S}^3$ vanishes identically (indeed, the equator is an isoperimetric surface of area $4\pi$), so that indeed $I(V) = I_{\mathbb{S}^3}(V)$ on $[0, V_0]$.
\end{proof}
\end{corollary}

In the next lemma we recall the Taylor series expansion for the volume of small geodesic balls, found by A. Gray and L. Vanhecke in \cite{GrayVanhecke} by integration of the expansion of the volume element in geodesic normal coordinates. For our purposes, the expansion of the volume element to fourth order (see also Lemma 3.4 in Chapter 5 in \cite{SchoenYau}) is sufficient.

\begin{lemma} [\cite{GrayVanhecke}] Let $(M^3, g)$ be a Riemannian $3$-manifold, $p \in M^3$, and $0 \leq r \ll
1$. Then \begin{equation*} \begin{split} &\vol(B(p, r)) = \\ &\frac {4 \pi r^3}{3} \Bigg (1 \underbrace {- \frac {\R(p)}{30}}_{:=c_1(p)} r^2 + \underbrace{\frac {1}{6300} \big(4 \R(p)^2 - 2|\Ric(p)|^2 - 9 (\Delta \R)(p) \big)}_{:=c_2(p)} r^4  + O(r^6)\Bigg).\end{split} \end{equation*}
\end{lemma}

In Section 8 of \cite{GrayVanhecke}, the authors concluded directly from this expansion that a $3$-manifold $M^3$ has constant sectional curvature (or, equivalently, is Einstein) provided that all small geodesic balls in $M^3$ have the same volume as the geodesic balls of the same radius in some fixed simply-connected $3$-manifold of constant sectional curvature. In our proof of Theorem \ref{thm: conclusion} we will use the isoperimetric profile $I_M$ as a \emph{lower} bound for the isoperimetric ratio of small geodesic balls, whose Taylor expansion we give in the lemma below. It turns out that this one-sided comparison (along with the lower bound for the scalar curvature) already implies that $M^3$ is round, and the proof of this fact given below depends delicately on the sign of the coefficients in the Taylor expansion of the preceding lemma.

\begin{lemma} Given $p \in M^3$ and $0 < V \ll 1$ there exists a unique $r > 0$
with $\vol (B(p, r)) = V$, and the dependence is smooth. Introducing the
variable $W = \big( \frac {3 V}{ 4\pi}\big)^{\frac 1 3}$ we obtain for the area
of the corresponding geodesic sphere the expansion
\begin{displaymath} \area(\partial B(p, r)) = 4 \pi W^2  \Big( 1 + c_1 W^2 +
\big(- \frac{11}{9} c_1^2 + \frac{5}{3} c_2 \big)W^4 + O(W^6) \Big),
\end{displaymath}
\noindent where $c_1 = c_1(p)$, $c_2 = c_2(p)$ are as in the previous lemma.
\begin{proof} This follows from the previous lemma and an elementary calculation. 
\end{proof}
\end{lemma}

\begin{theorem} \label{thm: conclusion} Let $(M^3, g)$ be a closed Riemannian
manifold with $\R_M \geq 6$ and $\Ric_M > 0$. Then all isoperimetric surfaces
of $M^3$ have area strictly less than $4\pi$ unless $(M^3, g)$ is isometric to
$\mathbb{S}^3$.

\begin{proof} Assume that $I(V_0) \geq 4\pi$ for some $V_0 \in (0, \vol(M^3)/2]$.
Then $m_H(V_0) \leq 0$ and Corollary \ref{cor: isoperimetric profiles agree} immediately implies that $\vol(M^3) = \vol(\mathbb{S}^3)$ and that $I_M = I_{\mathbb{S}^3}$.
\\

\noindent As a first step, we argue that $\R_M \equiv 6$. Fix a point $p
\in M^3$. For any (small) value of $V > 0$ the boundary of the geodesic ball
$B(p, r)$ of volume $V$ has surface area $\geq I_M(V) = I_{\mathbb{S}^3}(V)$.
Using that the isoperimetric surfaces of the sphere are just geodesic spheres,
we obtain from the previous lemma that \begin{equation*} 4 \pi W^2  \Big( 1 -
\frac {\R_M(p)}{30} W^2 + O(W^4) \Big) \geq 4 \pi W^2  \Big( 1 - \frac{1}{5}W^2
+ O(W^4) \Big), \end{equation*} where again $W = \big( \frac {3 V}{
4\pi}\big)^{\frac 1 3}$. Since this inequality holds for all $W > 0$
sufficiently small, we conclude that $\R_M(p) \leq 6$ and hence, since $p \in
M^3$ was arbitrary, that $\R_M \equiv 6$. \\

\noindent We now proceed to show that in fact $\Ric_M \equiv 2 g$. To see this,
we focus on the sixth order term in the above expansion: \begin {eqnarray*} 4
\pi W^2  \Big( 1 - \frac {1} {5} W^2 + \big( - \frac {|\Ric_M(p)|^2} {1890} -
\frac {17}{1575})W^4 + O(W^6) \Big) \\ \geq 4 \pi W^2  \Big( 1 - \frac {1} {5}
W^2 + \big( - \frac {|\Ric_{\mathbb{S}^3}(N)|^2} {1890} - \frac {17}{1575})W^4
+ O(W^6) \Big). \end {eqnarray*} Hence $|\Ric_M(p)|^2 \leq
|\Ric_{\mathbb{S}^3}(N)|^2 = 12$. Since $\R_M \equiv 6$ this means that $\Ric_M
\equiv 2g.$ \\

\noindent It follows that $M^3$ is round. Using now that $\vol(M^3) = \vol(\mathbb{S}^3)$ we see that $M^3$ must actually coincide with the standard
$\mathbb{S}^3$.
\end{proof}

\end{theorem}


\section{Min-Oo's Conjecture}

As an application of Theorem \ref{thm: conclusion} we answer in the affirmative the following special case of
Min-Oo's conjecture:

\begin{theorem} \label{thm: Min-Oo} Let $(M^3, g)$ be a compact orientable
Riemannian $3$-manifold with scalar curvature $\R_M \geq 6$, Ricci curvature $\Ric_M
> 0$ and totally geodesic boundary $\partial M^3$. If $\area (\partial M^3)
\geq 4 \pi$ and $\partial M^3$ is an isoperimetric surface of the doubled
manifold $(\hat M^3, \hat{g})$, then $(M^3, g)$ is isometric to the upper hemisphere
$\mathbb{S}^3_+$.

\begin{proof}
Observe that since we assume that $\partial M^3$ is totally geodesic in $M^3$, the
doubled manifold with the reflected metric is only $\mathcal{C}^{1, 1}$ across $\partial
M^3$. This degree of regularity is enough for the proof of Lemma \ref{lem:
monotonicity} to pertain. The full assertion now follows from applying the
method of Theorem \ref{thm: conclusion} away from the gluing region.
\end{proof}
\end{theorem}



\begin{thebibliography}{99}
\bibitem{AnderssonCaiGalloway} L. Andersson, M. Cai, G. Galloway: ``Rigidity
and positivity of mass for asymptotically hyperbolic manifolds," Ann. Henri
Poincar\'e 9 (2008). MR2389888
\bibitem{BavardPansu} C. Bavard, P. Pansu: ``Sur le volume minimal de $\mathbb{R}^2$,'' Ann. Sci. Ecole Norm. Sup. (4) 19 (1986). MR0875084 (88b:53048) 
\bibitem{HughBrayThesis} H. L. Bray: ``The Penrose inequality in general
relativity and volume comparison theorems involving scalar curvature," Thesis,
Stanford University (1997).  arXiv:0902.3241
\bibitem{ChristodoulouYau} D. Christodoulou, S.-T. Yau: ``Some remarks on the
quasi-local mass," Contemporary Mathematics, Volume 71 (1988). MR0954405 (89k:83050) 
\bibitem{ChruscielHerzlich} P.  Chrusciel, M. Herzlich: ``The mass of
asymptotically hyperbolic Riemannian manifolds," Pacific J. Math. 212 (2003). MR2038048 (2005d:53052) 
\bibitem{Corvino} J. Corvino: ``Scalar curvature deformation and a gluing
construction for the Einstein constraint equations,"  Comm. Math. Phys.  214
 (2000). MR1794269 (2002b:53050) 
\bibitem{GrayVanhecke} A. Gray, L. Vanhecke: ``Riemannian geometry as determined by the volumes of small geodesic balls,'' Acta. Math. 142 (1979). MR0521460 (81i:53038) 
\bibitem{HangWang} F. Hang, X. Wang: ``Rigidity and non-rigidity results on the
sphere,"  Comm. Anal. Geom.  14 (2006). MR2230571 (2007d:53065)
\bibitem{HangWang2} F. Hang, X. Wang: ``A rigidity theorem for the hemisphere,"
preprint (2007). arXiv:0711.4595v2
\bibitem{Huisken} G. Huisken: ``An isoperimetric concept for mass and quasilocal mass," Lecture at the \emph{Mathematical Aspects of General Relativity} workshop in Oberwolfach (2006), see http://www.mfo.de/programme/schedule/2006/02. 
\bibitem{Miao} P. Miao: ``Positive mass theorem on manifolds admitting corners
along hypersurfaces,'' Adv. Theor. Math. Phys. 6 (2002). MR1982695 (2005a:53065) 
\bibitem{Min-Oo} M. Min-Oo: ``Scalar curvature rigidity of asymptotically
hyperbolic spin manifolds,"  Math. Ann. 285  (1989). MR1027758 (91b:53047)
\bibitem{Min-Oo98} M. Min-Oo: ``Scalar curvature rigidity of certain symmetric
spaces," CRM Proceedings and Lecture Notes, Volume 15 (1998) MR1619128 (99e:53056) 
\bibitem{Ros} A. Ros: ``The isoperimetric problem," Lectures notes for the Clay
Mathematics Institute Summer School on the Global Theory of Minimal Surfaces, Amer.  Math. Soc. (2001). MR2167260 (2006e:53023)
\bibitem{SchoenYau} R. Schoen, S.-T. Yau: ``Lectures on Differential Geometry,"
Conference Proceedings and Lecture Notes in Geometry and Topology, Volume 1,
International Press (1994). MR1333601 (97d:53001)
\bibitem{PMTI} R. Schoen, S.-T. Yau: ``Positivity of the total mass of a
general spacetime,"  Phys. Rev. Lett.  43  (1979). MR0547753 (81c:58024) 
\bibitem{ShiTam} Y. Shi, L. Tam: ``Positive mass theorem and the boundary
behavior of compact manifolds with nonnegative scalar curvature," J.
Differential Geom. 62 (2002). MR1987378 (2005b:53046) 

\bibitem{Wang} X. Wang: ``The mass of asymptotically hyperbolic manifolds," J.
Differential Geom. 57 (2001). MR1879228 (2003c:53044) 
\bibitem{Witten} E. Witten: ``A new proof of the positive energy theorem,''
Comm. Math. Phys. 80 (1981). MR0626707 (83e:83035) 
\end{thebibliography}
\end{document}